
\documentstyle{amsppt}
\baselineskip18pt
\magnification=\magstep1
\pagewidth{30pc}
\pageheight{45pc}

\hyphenation{co-deter-min-ant co-deter-min-ants pa-ra-met-rised
pre-print pro-pa-gat-ing pro-pa-gate
fel-low-ship Cox-et-er dis-trib-ut-ive}
\def\leaderfill{\leaders\hbox to 1em{\hss.\hss}\hfill}

\def\tr{{\text {\rm \, tr}}}

\def\id{\text {\rm id}}

\def\sgn{{\text {\rm \, sgn}}}

\def\a{{\alpha}}
\def\be{{\beta}}

\def\e{{\varepsilon}}

\def\l{{\lambda}}

\def\bn{{\bold n}}

\def\bv{{\bold v}}

\def\b0{\text{\bf 0}}

\def\ra{{\ \longrightarrow \ }}

\def\hcpn{\text{\rm h}\gamma_n}

\def\lan{{\langle}}
\def\ran{{\rangle}}

\def\real{{\Bbb R}}
\def\complex{{\Bbb C}}
\def\zed{{\Bbb Z}}

\def\boxit#1{\vbox{\hrule\hbox{\vrule \kern3pt
\vbox{\kern3pt\hbox{#1}\kern3pt}\kern3pt\vrule}\hrule}}
\def\rabbit{\vbox{\hbox{\kern0pt
\vbox{\kern0pt{\hbox{---}}\kern3.5pt}}}}

\def\tableau#1{
        \hbox {
                \hskip -10pt plus0pt minus0pt
                \raise\baselineskip\hbox{
                \offinterlineskip
                \hbox{#1}}
                \hskip0.25em
        }
}

\def\tabCol#1{
\hbox{\vtop{\hrule
\halign{\strut\vrule\hskip0.5em##\hskip0.5em\hfill\vrule\cr\lower0pt
\hbox\bgroup$#1$\egroup \cr}
\hrule
} } \hskip -10.5pt plus0pt minus0pt}

\def\CR{
        $\egroup\cr
        \noalign{\hrule}
        \lower0pt\hbox\bgroup$
}



\def\blank#1#2{
\hbox to #1{\hfill \vbox to #2{\vfill}}
}


\def\strut{\vrule height10pt depth5pt width0pt}

\def\seca{1}
\def\secb{2}
\def\secc{3}
\def\secd{4}
\def\sym{{\frak S}}

\topmatter
\title Homology representations arising from the half cube, II
\endtitle

\author R.M. Green \endauthor
\affil Department of Mathematics \\ University of Colorado \\
Campus Box 395 \\ Boulder, CO  80309-0395 \\ USA \\ {\it  E-mail:}
rmg\@euclid.colorado.edu \\
\newline
\endaffil


\abstract
In a previous work, we defined a family of subcomplexes of the
$n$-dimensional half cube by removing the interiors of all half cube shaped
faces of dimension at least $k$, and we proved that the reduced 
homology of such a
subcomplex is concentrated in degree $k-1$.  This homology group supports a 
natural action of the Coxeter group $W(D_n)$ of type $D$.  In this paper,
we explicitly determine the characters (over ${\Bbb C}$) of these homology
representations, which turn out to be multiplicity free.  Regarded as
representations of the symmetric group $\sym_n$ by restriction, the homology
representations turn out to be direct sums of certain representations induced
from parabolic subgroups.  The latter representations of
$\sym_n$ agree (over ${\Bbb C}$) with the representations of 
$\sym_n$ on the $(k-2)$-nd homology of the complement of the $k$-equal 
real hyperplane arrangement.
\endabstract

\endtopmatter


\head \seca. Introduction \endhead

The half cube, also known as the demihypercube, 
may be constructed by selecting one point from each adjacent pair of 
vertices in the $n$-dimensional hypercube and taking the convex hull of 
the resulting set of $2^{n-1}$ points.  In a previous work \cite{{\bf 9}}, 
we showed that the $k$-faces of the half cube $\hcpn$ are of two types: 
regular simplices and, for $k \geq 3$, isometric copies of half cubes of 
lower dimension.  These faces assemble naturally into a regular CW complex,
$C_n$, which is homeomorphic to a ball.  Furthermore, for each 
$3 \leq k \leq n$, there is an interesting subcomplex $C_{n, k}$ of $C_n$
obtained by deleting the interiors of all the half cube shaped faces of 
dimensions 
$l \geq k$.  We also showed in \cite{{\bf 9}, Theorem 3.3.2} that the 
reduced homology of
$C_{n, k}$ is free over $\zed$ and concentrated in degree $k-1$.

The nonzero Betti numbers $B(n, k)$ of $C_{n, k}$ can be characterized by
simple recurrence relations: $B(n, 0) = B(n, n) = 1$ and, for 
$0 < k < n$, $$B(n, k) = 2 B(n-1, k) + B(n-1, k-1).
$$ There are also nonrecursive formulae for $B(n, k)$; for example, 
Bj\"orner--Welker \cite{{\bf 6}, Theorem 1.1 (c)} prove that $$
B(n, k) = \sum_{i = k}^n {n \choose i} {{i-1} \choose {k-1}}
,$$ where we interpret ${{-1} \choose {-1}}$ to mean $1$.
The numbers $B(n, k)$ are interesting because they occur
in a diverse range of contexts, such as:
\item{(i)}{in the problem of finding, given $n$ real numbers,
a lower bound for the complexity 
of determining whether some $k$ of them are equal 
\cite{{\bf 4}, {\bf 5}, {\bf 6}, \S1},}
\item{(ii)}{as the $(k-2)$-nd Betti numbers of the $k$-equal real hyperplane
arrangement in $\real^n$ \cite{{\bf 6}},}
\item{(iii)}{as the ranks of $A$-groups appearing in combinatorial homotopy
theory \cite{{\bf 1}, {\bf 2}},}
\item{(iv)}{as the number of nodes used by the Kronrod--Patterson--Smolyak
cubature formula in numerical analysis \cite{{\bf 15}, Table 3}, and}
\item{(v)}{(when $k = 3$) in engineering, as the number of three-dimensional 
block structures associated to $n$ joint systems in the construction
of stable underground structures \cite{{\bf 12}}.}

The connections between (i), (ii) and (iii) above are now well understood.
Although the half cube polytope has no obvious direct relationship with any
of these five phenomena, its associated homology groups share an intriguing 
feature in common with those appearing in (ii) and (iii): they all
support natural actions of the symmetric group $\sym_n$.  One possible way 
to forge a link between the half cube and the situations in (ii) or (iii) 
is to try to understand the various homology groups in terms of group 
representations.

The $k$-equal real hyperplane arrangement 
$V_{n, k}^\real$ is the set of points $(x_1, \ldots, x_n) \in \real^n$
such that $x_{i_1} = x_{i_2} = \cdots = x_{i_k}$ for some set of indices
$1 \leq i_1 < i_2 < \cdots < i_k \leq n$.  The 
complement $\real^n - V_{n, k}^\real$, denoted by $M_{n, k}^\real$, is
a manifold whose homology is concentrated in degrees $t(k-2)$, where 
$t \in \zed$ satisfies $0 \leq t \leq {\lfloor {n \over k} \rfloor}$ 
(see \cite{{\bf 6}, Theorem 1.1(b)}).  It is clear from the definition that
$M_{n, k}^\real$ supports an action of $\sym_n$ via permutation of coordinates,
and this action endows the nonzero homology groups with the structure of
$\sym_n$-modules.  The characters of these modules were computed explicitly 
by Peeva, Reiner and Welker in \cite{{\bf 14}, Theorem 4.4}.

The $n$-dimensional half cube has a large symmetry group $G_n$ of orthogonal
transformations acting on it via cellular automorphisms.  This group always 
contains the Coxeter group $W(D_n)$ of order $2^{n-1}n!$, although this 
containment is proper if $n = 4$, and in turn the group $W(D_n)$ contains
a subgroup isomorphic to $\sym_n$.  The action of $W(D_n)$ on $\hcpn$ 
induces, for each $k$, an action on the nonzero homology groups of 
$C_{n, k}$.  In this paper, we will compute the characters (over $\complex$) 
of these homology representations.  Regarded as modules for $W(D_n)$, the
representations turn out to be multiplicity free (Theorem \secd.4), although
they are not generally induced modules in any nontrivial sense.  In contrast,
if the homology representations are regarded as modules for the symmetric 
group $\sym_n$ by restriction, then they are no longer multiplicity free,
but they do turn out to be isomorphic to direct sums of modules induced from 
maximal Young subgroups (Theorem \secd.7).  
Furthermore, over the complex numbers, the action of $\sym_n$ on the
$(k-1)$-st homology of $C_{n, k}$ agrees with the action of $\sym_n$ on 
the $(k-2)$-nd homology of $M_{n, k}^\real$.

Our results have some interesting combinatorial consequences.  One of
these (Corollary \secd.6) is that if we restrict the
representation of $W(D_n)$ on the $(k-1)$-st homology of $C_{n, k}$ 
to the subgroup $W(D_{n-1})$, then the corresponding branching
rule categorifies the usual recurrence relation for the Betti numbers 
$B(n, k)$.  Another nice property is that if one computes the
dimension of the homology representations from a knowledge of their 
characters, then one obtains a combinatorial proof of the Bj\"orner--Welker
formula for $B(n, k)$ mentioned above.

\head \secb. Character theory of Coxeter groups of classical type \endhead

The main groups of interest in this paper are the finite Coxeter groups of
classical type, meaning types $A_{n-1}$, $B_n$ and $D_n$.  It will be 
convenient to number the vertices of the corresponding Coxeter graphs as 
shown in Figure 1.

\topcaption{Figure 1} Coxeter graphs of type $A_{n-1}$, $B_n$ and $D_n$
\endcaption
\centerline{
\hbox to 3.458in{
\vbox to 1.638in{\vfill
        \includegraphics{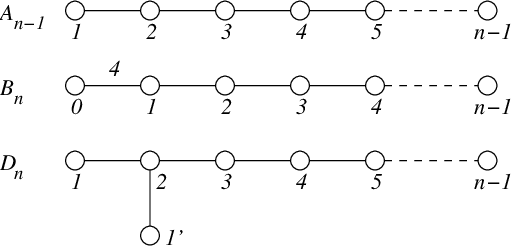}
}
\hfill}
}

We now summarize some well-known properties of these groups.  More details
may be found in \cite{{\bf 11}} or \cite{{\bf 3}}.

The {\it Coxeter group} $W = W(\Gamma)$ corresponding to a Coxeter graph 
$\Gamma$ with vertices $S = S(\Gamma)$ is given by the presentation $$
\langle s_i : i \in S(\Gamma) : (s_i s_j)^{m_{ij}} = 1 \rangle.
$$ The numbers $m_{ij}$ are defined to satisfy $m_{ii} = 1$ and 
$m_{ij} = m_{ji}$ for all $i, j \in S$.  Furthermore, we have $m_{ij} = 2$
if $i$ and $j$ are not adjacent in the graph; $m_{ij} = 3$ if $i$ and $j$
are connected by an unlabelled edge; and $m_{ij} = k$ if $i$ and $j$ are
connected by an edge labelled $k > 3$.  

If $S' \subset S$, then we refer to the subgroup $W'$ of $W$ that is 
generated by $S'$ as a {\it parabolic subgroup} of $W$.  In this case, 
$W'$ inherits the structure of a Coxeter group from $W$.

The Coxeter group $W(A_{n-1})$ is isomorphic (as an abstract group) to the
symmetric group $\sym_n$, and the Coxeter generator $s_i$ may be identified
with the transposition $(i, i+1)$.

The Coxeter group $W(B_n)$ is isomorphic to the wreath product 
$\zed_2 \wr \sym_n$ of order $2^n n!$.  This may be regarded as a group of
permutations of $n$ signed objects, in which $s_i$ acts by the transposition
$(i, i+1)$ for $1 \leq i < n$, and $s_0$ acts by changing the sign of the
object numbered $1$.  The parabolic subgroup of $W(B_n)$ obtained by omitting
the generator $s_0$ is canonically isomorphic (as a Coxeter group) to 
$W(A_{n-1})$.

The Coxeter group $W(D_n)$, which will be our main group of interest, 
can be identified with the index $2$ subgroup of 
$W(B_n)$ consisting of those elements effecting an even number of sign
changes.  As before, we may identify $s_i$ with the permutation $(i, i+1)$
for $1 \leq i < n$.  The other generator, $s_{1'}$, can be identified with
the element $s_0 s_1 s_0$ of $W(B_n)$.  It therefore acts by changing the
sign of each of objects $1$ and $2$, followed by the transposition $(1, 2)$.
The parabolic subgroup of $W(D_n)$ obtained by omitting the generator $s_{1'}$
is canonically isomorphic (as a Coxeter group) to $W(A_{n-1})$.  We will
often abuse notation slightly and refer to this subgroup as $\sym_n$.

For any Coxeter group $(W, S)$, there is a unique homomorphism 
$\e : W \ra \{\pm 1\}$ to the multiplicative group of $2$ elements sending 
each element of $S$ to $-1$.  This homomorphism is known as the {\it sign
representation}.  We will write $\sgn_n$ (respectively, $\id_n$) to denote 
the sign (respectively, trivial) representation of
any of the Coxeter groups of types $A_{n-1}$, $B_n$ or $D_n$.

The character theory of finite Coxeter groups of classical type is well 
understood.  It is described in Geck and Pfeiffer's book \cite{{\bf 7}, \S5} 
and, more explicitly, in Stembridge's notes on the topic \cite{{\bf 17}}.
We now summarize some of the key properties of the theory for later use.

The irreducible representations of $W(A_{n-1})$ (over $\complex$) are 
indexed by the partitions of $n$, or equivalently, the set of Young diagrams
of size $n$.  We will write the corresponding set of characters as $$
\{\chi^\l : |\l| = n\}
.$$  The degree, $\chi^\l(1)$, of $\chi^\l$ is the number of standard 
Young tableaux
of shape $\l$; that is, the number of ways of filling a Young diagram of
shape $\l$ with the numbers $1, 2, \ldots, n$ once each in such a way that
the entries increase along rows and down columns.  The identity character
corresponds to the partition $\l = [n]$, whose Young diagram has one row, 
and the sign character corresponds to the partition $\l = [1^n]$, whose Young
diagram has one column.  

Another important character for our purposes
corresponds to the partition $$[n-1, 1],$$ which gives the character of 
the {\it reflection representation} associated to the Coxeter group of type 
$A_{n-1}$.  This representation may also be constructed by first
taking the $n$-dimensional representation of $\sym_n$ corresponding to the 
natural action of the group on $n$ letters, and then quotienting by the 
$1$-dimensional submodule spanned by the all-ones vector.  

The irreducible characters of $W(B_n)$ are indexed by the set $$
\{\chi^{(\mu, \nu)} : |\mu| + |\nu| = n\}
.$$  The dimensions of the irreducibles may be obtained from the corresponding
dimensions in type $A$ via the formula $$
\chi^{(\mu, \nu)}(1) = {n \choose {|\mu|}} \chi^\mu(1) \chi^\nu(1)
.$$  The identity character corresponds to the pair $([n], [0])$, and the sign
character to the pair $([0], [1^n])$.

As described above, we may regard $W(D_n)$ as a subgroup of $W(B_n)$.  Under
this identification, the irreducible characters $\chi^{(\mu, \nu)}$ and
$\chi^{(\nu, \mu)}$ (where $\mu \ne \nu$) both restrict to the same 
irreducible character of $W(D_n)$, which we denote by $\chi^{\{\mu, \nu\}}$.
On the other hand, the irreducible character $\chi^{(\mu, \mu)}$ of $W(B_n)$
restricts to a sum of two nonisomorphic irreducible characters of
$W(D_n)$ of the same degree; we denote the latter by $\chi_+^{\{\mu, \mu\}}$ 
and $\chi_-^{\{\mu, \mu\}}$.  These exhaust all the irreducible characters of
$W(D_n)$. In other words, 
the irreducible characters of $W(D_n)$ are indexed by the set $$
\{\chi^{\{\mu, \nu\}} : |\mu| + |\nu| = n\} \cup \{\chi_\pm^{\{\mu, \mu\}}
: |\mu| = n/2\}
,$$ where the second subset is empty if $n$ is odd.  The identity character
corresponds to the pair $\{[n], [0]\}$ and the sign character corresponds 
to the pair $\{[1^n], [0]\}$.  It is immediate from
the above remarks that the dimensions of the corresponding irreducibles are 
given by $$
\chi^{\{\mu, \nu\}}(1) = {n \choose {|\mu|}} \chi^\mu(1) \chi^\nu(1)
$$ and $$
\chi_\pm^{\{\mu, \mu\}}(1) = {1 \over 2} {n \choose {|\mu|}} \chi^\mu(1)^2
.$$

The following two lemmas concerning characters of $W(D_n)$ will be important
in the sequel.  It will sometimes be convenient to write 
$\chi_\e^{\{\mu, \nu\}}$ to refer to the irreducible character 
$\chi^{\{\mu, \nu\}}$ if $\mu \ne \nu$, and to refer to either of the 
irreducible characters $\chi_+^{\{\mu, \nu\}}$ or 
$\chi_-^{\{\mu, \nu\}}$ if $\mu = \nu$.

\proclaim{Lemma \secb.1}
Let $\sym_n$ be the parabolic subgroup of $G = W(D_n)$ obtained by omitting 
the generator $s_{1'}$.  
Let $m = {\lfloor{n \over 2} \rfloor}$.
\item{\rm (i)}{If $\mu \ne \nu$, then we have $$
\chi^{\{\mu, \nu\}} \downarrow^G_{\sym_n} = \sum_\lambda c_{\mu \nu}^\lambda
\chi^\lambda
,$$ where the $c_{\mu \nu}^\lambda$ are the Littlewood--Richardson 
coefficients.}
\item{\rm (ii)}{If $n$ is odd, then we have $$
\id_n \uparrow_{\sym_n}^G = \sum_{l \leq m} \chi^{\{[l], [n - l]\}}
.$$}
\item{\rm (iii)}{If $n$ is even, then we have $$
\id_n \uparrow_{\sym_n}^G = 
\chi_+^{\{[m], [m]\}} +
\sum_{l < m} \chi^{\{[l], [n - l]\}}
.$$}
\endproclaim

\demo{Proof}
Part (i) appears in \cite{{\bf 17}, \S3A}.
Under the hypotheses of part (ii), we must have $\mu \ne \nu$ because $n$ is
odd.  The conclusion of (ii) then follows from (i) by using Frobenius
reciprocity and the Pieri Rule.
Part (iii) appears in \cite{{\bf 17}, \S3C}.
\qed\enddemo

\proclaim{Lemma \secb.2}
Let $G = W(D_n)$, let $3 \leq k \leq n$, and let $D_k$ (respectively, 
$D_{n-k}$) be the parabolic subgroup of $G$ generated by the set $$
\{s_{1'}\} \cup \{s_i : 1 \leq i < k\}
$$ (respectively, $\{s_i : i > k\}$.)  Denote the usual inner product on
characters by $\lan , \ran$.  Suppose below that the ordered pairs $(\a, \psi)$
and $(\be, \theta)$ are not equal.
\item{\rm (i)}{If $\mu \ne \nu$ then we have $$
\left\langle
\chi^{\{\mu, \nu\}} \downarrow^G_{D_k \times D_{n-k}}, \ 
\chi_\e^{\{\a, \be\}} 
\times
\chi_{\e'}^{\{\psi, \theta\}} 
\right\rangle = 
c_{\a \psi}^\mu
c_{\be \theta}^\nu
+
c_{\a \theta}^\mu
c_{\be \psi}^\nu
+ 
c_{\be \psi}^\mu
c_{\a \theta}^\nu
+
c_{\be \theta}^\mu
c_{\a \psi}^\nu
.$$}
\item{\rm (ii)}{We have $$
\left\langle
\chi_{\pm}^{\{\mu, \mu\}} \downarrow^G_{D_k \times D_{n-k}}, \ 
\chi_\e^{\{\a, \be\}} 
\times
\chi_{\e'}^{\{\psi, \theta\}} 
\right\rangle = 
c_{\a \psi}^\mu
c_{\be \theta}^\mu
+
c_{\a \theta}^\mu
c_{\be \psi}^\mu
.$$}
\endproclaim

\demo{Proof}
After applying Frobenius reciprocity, this becomes a restatement of results 
in \cite{{\bf 17}, \S3A, \S3D}.
\qed\enddemo

The irreducible characters in type $B_n$ have the following well-known 
branching rule, which will be useful in the sequel.

\proclaim{Lemma \secb.3}
Let $\chi^{(\lambda, \mu)}$ be an irreducible character for $W(B_n)$.  Then
we have $$
\chi^{(\lambda, \mu)} \downarrow^{W(B_n)}_{W(B_{n-1})}
= 
\sum_{d \in I(\lambda)} 
\chi^{(\lambda^{(d)}, \mu)}
+ 
\sum_{d \in I(\mu)} 
\chi^{(\lambda, \mu^{(d)})}
,$$ where $I(\lambda)$ is the set of removable boxes in the Young diagram
corresponding to $\lambda$, and $\lambda^{(d)}$ is the result of removing
box $d$ from the Young diagram.
\endproclaim

\demo{Proof}
This appears in \cite{{\bf 7}, \S6.1.9}.
\qed\enddemo

\proclaim{Corollary \secb.4}
Maintain the notation of Lemma \secb.3.  Suppose in addition that each 
character $\chi^{(\a, \be)}$ appearing in Lemma \secb.3 satisfies 
$\a \ne \be$.  Then we have $$
\chi^{\{\lambda, \mu\}} \downarrow^{W(D_n)}_{W(D_{n-1})}
= 
\sum_{d \in I(\lambda)} 
\chi^{\{\lambda^{(d)}, \mu\}}
+ 
\sum_{d \in I(\mu)} 
\chi^{\{\lambda, \mu^{(d)}\}}
.$$
\endproclaim

\demo{Proof}
Recall that each type $B$ character $\chi^{(\a, \be)}$ appearing in the 
statement restricts to the irreducible type $D$ character 
$\chi^{\{\a, \be\}}$.
The result now follows from Lemma \secb.3 and the fact that, under the
usual identifications, we have $W(B_{n-1}) \cap W(D_n) = W(D_{n-1})$.
\qed\enddemo

\head \secc. The half cube \endhead

An $n$-dimensional {\it (Euclidean) polytope} 
$\Pi_n$ is a closed, bounded, convex subset of $\real^n$ enclosed by a 
finite number of hyperplanes.  The part of $\Pi_n$ that lies in one of
the hyperplanes is called a {\it facet}, and each facet is an 
$(n-1)$-dimensional polytope.  
A polytope is homeomorphic to an $n$-ball (which follows, for example, from
\cite{{\bf 13}, Lemma 1.1}), and the boundary of the polytope,
which is equal to the union of its facets, is identified with the 
$(n-1)$-sphere by this homeomorphism.

Iterating this construction gives
rise to a set of $k$-dimensional polytopes $\Pi_k$ (called {\it $k$-faces}) 
for each $0 \leq k \leq n$.  The elements of $\Pi_0$ are called {\it vertices}
and the elements of $\Pi_1$ are called {\it edges}.  It is not hard to show
that a polytope is the convex hull of its set of vertices, and that the 
boundary of a polytope is precisely the union of its $k$-faces for 
$0 \leq k < n$.  What is less obvious, but still true 
\cite{{\bf 18}, Theorem 1.1}, is that the convex hull of an arbitrary 
finite subset
of ${\Bbb R}^n$ is a polytope in the above sense.  It follows that a polytope
is determined by its vertex set, and we write $\Pi(V)$ for the polytope 
whose vertex set is $V$.

\definition{Definition \secc.1}
Let $n \geq 4$ be an integer, and let $\bn = \{1, 2, \ldots, n\}$.

We define the set $\Psi_n$ to be the set of
$2^n$ vertices of the hypercube whose coordinates are $$
(\pm 1, \pm 1, \ldots, \pm 1)
.$$ Let $\Psi^+_n$ be the set of
$2^{n-1}$ vertices with an even number of negative coordinates, 
and let $\Psi^-_n$ be $\Psi_n \backslash \Psi^+_n$.

Let $\bv' \in \Psi^-_n$ and $S \subseteq \bn$.  We define the subset 
$K(\bv', S)$
of $\Psi^+_n$ by the condition that $\bv \in K(\bv', S)$ if and only if 
$\bv$ and $\bv'$ differ only in the $i$-th coordinate, and $i \in S$.

Let $\bv \in \Psi^+_n$ and let $S \subseteq \bn$.
We define the subset $L(\bv, S)$ of $\Psi^+_n$ by the 
condition that $\bv' \in L(\bv, S)$ if and only if $\bv$ and $\bv'$ agree
in the $i$-th coordinate whenever $i \not\in S$.  The set $S$ is 
characterized as the set of coordinates
at which not all points of $L(\bv, S)$ agree.
\enddefinition

The $k$-faces of the half cube were classified in \cite{{\bf 9}}.

\proclaim{Theorem \secc.2 (\cite{{\bf 9}})}
The $k$-faces of $\hcpn$ for $k \leq n$ are as follows:
\item{\rm (i)}{$2^{n-1}$ $0$-faces (vertices) given by the elements of
$\Psi^+_n$;}
\item{\rm (ii)}{$2^{n-2} {n \choose 2}$ $1$-faces $\Pi(K(\bv', S))$,
where $\bv' \in \Psi^-_n$ and $|S| = 2$;}
\item{\rm (iii)}{$2^{n-1} {n \choose 3}$ simplex shaped $2$-faces 
$\Pi(K(\bv', S))$, where $\bv' \in \Psi^-_n$ and $|S| = 3$;}
\item{\rm (iv)}{$2^{n-1} {n \choose {k+1}}$ simplex shaped $k$-faces 
$\Pi(K(\bv', S))$, where $\bv' \in \Psi^-_n$ and $|S| = k+1$ for 
$3 \leq k < n$;}
\item{\rm (v)}{$2^{n-k} {n \choose k}$ half cube shaped $k$-faces 
$\Pi(L(\bv, S))$, where $\bv \in \Psi^+_n$ and $|S| = k$ for 
$3 \leq k \leq n$.}

Furthermore, two faces are conjugate under the action of $W(D_n)$ if and
only if they have the same dimension and the same shape.
\endproclaim

\demo{Proof}
The classification of the $k$-faces is given in \cite{{\bf 9}, Theorem 2.3.6},
and the classification of the orbits under the action of $W(D_n)$ is given
in \cite{{\bf 9}, Theorem 4.2.3 (ii)}.
\qed\enddemo

The unique $n$-face in (v) above corresponds to the interior of the polytope.
The $k$-faces assemble naturally into a regular CW complex, $C_n$.  

\definition{Definition \secc.3}
For each integer $k$ with $-1 \leq k \leq n$, let $V_k$ be the $k$-th 
chain group in the complex $C_n$.  For $k \geq 3$,
we write $V_k = X_k \oplus Y_k$, where $X_k$ (respectively, $Y_k$) is the 
span of the simplex-shaped (respectively, half-cube shaped) faces.  If
$-1 \leq k < 3$, we define $X_k = V_k$ and $Y_k = 0$.
\enddefinition

We now recall some of the key properties of this complex; the reader is 
referred to \cite{{\bf 9}} for full details.  

For any fixed $k$ such that $3 \leq k \leq n$, one may form a CW 
subcomplex $C_{n, k}$ by removing the interiors of all the half cube 
shaped $l$-faces for $l \geq k$.  In other words, the $l$-th chain group 
of $C_{n, k}$ is equal to $V_k$ if $l < k$, and to $X_k$ if $l \geq k$.
The reduced (cellular) homology of 
$C_{n, k}$ is free over $\zed$ and concentrated in degree $k-1$ 
\cite{{\bf 9}, Theorem 3.3.2}.

The Coxeter group $W(D_n)$ acts naturally on $\Psi^+_n$ 
(and also on $\Psi^-_n$)
via signed permutations of the coordinates.  This induces an action of
$W(D_n)$ on the half cube $\hcpn$ via cellular automorphisms.  In particular,
elements of $W(D_n)$ send $k$-faces of $\hcpn$ to other $k$-faces of the same
type (i.e., simplex shaped or half cube shaped), which means that the
free $\zed$-modules $X_k$ and $Y_k$ of Definition \secc.3 acquire the
structure of $W(D_n)$-modules.
In turn, there is an induced action of $W(D_n)$ on the subcomplex 
$C_{n, k}$ via cellular automorphisms, as well as on the homology groups 
of $C_{n, k}$ \cite{{\bf 9}, Theorem 4.2.3}.

The following basic result will be of key importance in the sequel.

\proclaim{Lemma \secc.4}
Let $n \geq 3$ and let $s$ be a Coxeter generator of the group $G = W(D_n)$. 
The element $s$ acts on the half cube $\hcpn$ by a reflection in a
hyperplane through the origin.  The induced action of $s$ on 
$H_{n-1}(C_{n, n})$ and on the $n$-th chain group of $C_n$ is negation.
\endproclaim

\demo{Proof}
The first assertion follows from \cite{{\bf 8}, Proposition 3.6, Lemma 5.3}.

The CW space corresponding to the subcomplex $C_{n, n}$ is obtained from 
$C_n$ by deleting the (interior of the) unique $n$-cell.  It follows that
this space is homeomorphic to $S^{n-1}$.  A well-known result
\cite{{\bf 10}, 2.2 (e)} then shows that $s$ acts on $H_{n-1}(S^{n-1}; \zed)$
by negation, proving the first part of the second assertion.  For the final
assertion, we use
the fact that $s$ acts continuously on $\hcpn$, which means that it acts
on the chain complex of $C_n$ by a chain map.  Since the $n$-th chain group
of $C_n$ has rank $1$, the fact that $s$ acts by negation on 
$H_{n-1}(S^{n-1}; \zed)$ forces it to act by negation on $C_n$,
which completes the proof.
\qed\enddemo

\proclaim{Lemma \secc.5}
Let $n \geq 4$, $G = W(D_n)$ and let $k$ satisfy $3 \leq k \leq n$.
Let $D_k$ denote the parabolic subgroup of $W(D_n)$ generated by the set $$
\{s_{1'}\} \cup \{s_i : 1 \leq i < k\}
,$$ and let $\sym_{n-k}$ denote the parabolic subgroup of $W(D_n)$ generated
by the set $$
\{s_i : i > k\}
.$$  Regarding $Y_k$ as a $\complex G$-module by extension of scalars, we 
have $$
Y_k \cong_G (\sgn_k \otimes \id_{n-k})\uparrow_{D_k \times \sym_{n-k}}^G
.$$
\endproclaim

\demo{Proof}
By Theorem \secc.2, there is one orbit of half cube shaped faces for
each $3 \leq k \leq n$.  One of these has vertex set $L(\bv, S)$, where $$
\bv = (1, 1, \ldots, 1)
$$ and $$
S = \{k+1, k+2, \ldots, n\}
.$$  Let $e$ be the $k$-cell of the CW complex $C_n$ corresponding to 
$L(\bv, S)$.  It is clear from the definitions of the action of $G$ as signed
permutations that the set $L(\bv, S)$ is fixed setwise by all the $s_i$ other
than $s_k$.  The group generated by this subset of the generators is
$G_k := D_k \times \sym_{n-k}$, which has order $2^{k-1} k! (n-k)!$ and 
index $$
2^{n-k} {n \choose k}
$$ in $G$.  It now follows from Theorem \secc.2 (v) that $G_k$
is the full set stabilizer of $L(\bv, S)$.

The Coxeter generators $s_{1'}$, $s_1, \ldots, s_{k-1}$ of $D_k$ act
as reflections in hyperplanes through the origin.  Lemma \secc.4 then
shows that each of these generators sends $e$ to $-e$.
In contrast, the Coxeter generators $s_{k+1}, s_{k+2}, \ldots, s_n$ fix
$L(\bv, S)$ (and its convex hull) pointwise.  These generators fix $e$.

The assertion follows from the above observations.
\qed\enddemo

\proclaim{Lemma \secc.6}
Maintain the notation of Lemma \secc.5, and let $\eta(k, e)$ denote 
the partition $[e+1, 1, \ldots, 1]$ of $k+e$.
The character of the module $Y_k$ is given by $$
\sum_{e \leq n-k} \chi^{\{\eta(k, e), [n-k-e]\}} + 
\sum_{e' \leq n-(k+1)} \chi^{\{\eta(k+1, e'), [n-(k+1)-e']\}}
.$$
\endproclaim

\demo{Proof}
By transitivity of induction and Lemma \secc.5, we have $$
Y_k \cong_G 
\left( (\sgn_k \otimes \id_{n-k})
\uparrow_{D_k \times \sym_{n-k}}^{D_k \times D_{n-k}} \right)
\uparrow_{D_k \times D_{n-k}}^G
.$$

Let $m = {\lfloor{{n-k} \over 2} \rfloor}$.
By Lemma \secb.1, we have $$
(\sgn_k \otimes \id_{n-k})
\uparrow_{D_k \times \sym_{n-k}}^{D_k \times D_{n-k}}
= \sum_{l \leq m} \chi^{\{[1^k], [0]\}} \times \chi^{\{[l], [n-k-l]\}}
$$ if $n-k$ is odd, and $$
(\sgn_k \otimes \id_{n-k})
\uparrow_{D_k \times \sym_{n-k}}^{D_k \times D_{n-k}}
= 
\chi^{\{[1^k], [0]\}} \times \chi_+^{\{[m], [m]\}}
+ 
\sum_{l < m} \chi^{\{[1^k], [0]\}} \times \chi^{\{[l], [n-k-l]\}}
$$ if $n-k$ is even.  Lemma \secb.2 is applicable in this situation,
because neither of the partitions $[1^k]$ or $[0]$ has one row.
The assertion now follows from Lemma \secb.2 and the Pieri rule.
(Observe that the numbers appearing in Lemma \secb.2 (ii) are always zero
in this case.)
\qed\enddemo

\remark{Remark \secc.7}
The methods used in Lemma \secc.6 to determine the
characters of the modules $Y_k$ can be extended to compute the characters
of the modules $X_k$, as well as the characters of all the representations
corresponding to cycles and to boundaries in the subcomplexes $C_{n, k}$.
\endremark

\head \secd. Main results \endhead

In order to prove our main results, we require a version of the Hopf trace
formula that applies in contexts more general than simplicial complexes.

\proclaim{Theorem \secd.1 (Hopf trace formula \cite{{\bf 13}, Theorem 22.1})}
Let $K$ be a finite complex with (integral) chain groups $C_p(K)$ and
homology groups $H_p(K)$.  Let $T_p(K)$ be the torsion subgroup of
$H_p(K)$.  Let $\phi : C_p(K) \ra C_p(K)$ be a chain map, and let
$\phi_*$ be the induced map on homology.  Then we have $$
\sum_p (-1)^p \tr(\phi, C_p(K)) = 
\sum_p (-1)^p \tr(\phi_*, H_p(K)/T_p(K))
.$$ \qed
\endproclaim

\proclaim{Lemma \secd.2}
Consider the CW complex $C_n$ of Definition 
\secc.3; its chain groups are the $V_l$ for $-1 \leq l \leq n$.  
Let $\phi$ be a chain map of this chain complex.  Then we have $$
\sum_p (-1)^p \tr(\phi, V_p) = 0
.$$
\endproclaim

\demo{Proof}
The chain complex $C_n$ is a CW decomposition of the half cube, which is
a contractible space and has trivial reduced homology.  Theorem \secd.1
applies to the complex $C_n$, and the previous observation shows that the 
right hand side of the Hopf trace formula is zero, completing the proof.
\qed\enddemo

\proclaim{Lemma \secd.3}
Consider the CW subcomplex $C_{n, k}$ of
$C_n$; its chain groups are $V_l$ for $-1 \leq l < k$ and $X_l$ for
$k \leq l \leq n$, where $V_l = X_l \oplus Y_l$ for $l \geq 3$.
Let $\phi$ be a chain map of this chain complex.  Then we have $$
\tr(\phi_*, H_{k-1}(C_{n, k})) = 
\sum_{l \geq k} (-1)^{l-k} \tr(\phi, Y_l)
.$$
\endproclaim

\demo{Proof}
We first apply the Hopf trace formula to $C_{n, k}$ to obtain $$
\sum_p (-1)^p \tr(\phi, C_p(C_{n, k})) = 
\sum_p (-1)^p \tr(\phi_*, H_p(C_{n, k}))
;$$ there is no torsion because the homology of $C_{n, k}$ is free over
$\zed$ by \cite{{\bf 9}, Theorem 3.3.2}.  Since, by the same result, 
the reduced homology of $C_{n, k}$ is concentrated in degree
$k-1$, this simplifies to $$
(-1)^{k-1} \tr(\phi_*, H_{k-1}(C_{n, k})) = 
\sum_p (-1)^p \tr(\phi, C_p(C_{n, k}))
.$$  By Lemma \secd.2, we have $$
\sum_p (-1)^p \tr(\phi, C_p(C_{n, k})) +
\sum_{p \geq k} (-1)^p \tr(\phi, Y_k) = 0
,$$ which, combined with the preceding equation, gives $$
(-1)^{k-1} \tr(\phi_*, H_{k-1}(C_{n, k})) = 
\sum_{p \geq k} (-1)^{p+1} \tr(\phi, Y_k) = 0
.$$  The assertion now follows by multiplying both sides by $(-1)^{k-1}$.
\qed\enddemo

\proclaim{Theorem \secd.4}
Let $n \geq 4$, $G = W(D_n)$ and let $k$ satisfy $3 \leq k \leq n$.
Let $\eta(k, e)$ denote the partition $[e+1, 1, \ldots, 1]$ of $k+e$.
The character of the representation of $G$ on the $(k-1)$-st homology
of the complex $C_{n, k}$ is given by $$
\chi_D(n, k) = \sum_{e \leq n-k} \chi^{\{\eta(k, e), [n-k-e]\}}
.$$
\endproclaim

\demo{Proof}
Let $\chi_k$ denote the character of the $G$-module $Y_k$.  By Lemma 
\secd.3, the character of the homology representation is given by the
alternating sum $$
\chi_k - \chi_{k+1} + \chi_{k+2} - \chi_{k+3} \cdots
.$$  The result now follows from Lemma \secc.6: all of the terms appearing
in the statement of that result cancel, except those involving a partition
of the form $\eta(l, e)$ for $l = k$.
\qed\enddemo

\remark{Remark \secd.5}
It is known \cite{{\bf 16}, \S4} that the $k$-th exterior
power of the ($n$-dimensional) reflection representation of $W(D_n)$ is
irreducible and corresponds to the pair of partitions $$
\{ [1^k], [n-k] \} = \{\eta(k, 0), [n-k-0]\}
.$$  Theorem \secd.4 shows that this is one of the constituents of the
representation of $W(D_n)$ on the $(k-1)$-st homology of $C_{n, k}$.
\endremark

\proclaim{Corollary \secd.6}
Maintain the notation of Theorem \secd.4, and assume that $k < n$.  We have $$
\chi_D(n, k) 
\downarrow^{W(D_n)}_{W(D_{n-1})} = 
2 \chi_D(n-1, k) + 
\chi_D(n-1, k-1) 
.$$
\endproclaim

\demo{Proof}
Since $k \geq 3$, Corollary \secb.4 shows that $$\eqalign{
\chi^{\{\eta(k, e), [n-k-e]\}} \downarrow^{W(D_n)}_{W(D_{n-1})}
=& 
\chi^{\{\eta(k, e-1), [n-k-e]\}}
+ 
\chi^{\{\eta(k, e), [n-k-e-1]\}}\cr
&+ 
\chi^{\{\eta(k-1, e), [n-k-e]\}}\cr
=&
\chi^{\{\eta(k, e-1), [(n-1)-k-(e-1)]\}}
+ 
\chi^{\{\eta(k, e), [(n-1)-k-e]\}}\cr
&+ 
\chi^{\{\eta(k-1, e), [(n-1)-(k-1)-e]\}}, \cr
}$$ where we ignore any terms involving partitions with negative parts.
The result now follows by summing over $e$, as in Theorem \secd.4.
\qed\enddemo

\proclaim{Theorem \secd.7}
Let $n \geq 4$ and let $k$ satisfy $3 \leq k \leq n$.
Let $\sym_n$ denote the parabolic subgroup of $W(D_n)$ corresponding to the
omission of the generator $s_{1'}$.  Let $E_l$ be the 
$(l-1)$-dimensional reflection representation for $\sym_l$ described 
in \S\secb.  
\item{\rm (i)}{Regarded as a $\complex\sym_n$-module
by restriction, the $(k-1)$-st homology of the complex $C_{n, k}$ is 
isomorphic to $$
\bigoplus_{e \leq n-k}
\big( 
\id_{n-k-e} \otimes \bigwedge^{k-1} E_{k+e}
\big) \uparrow_{\sym_{n-k-e} \times \sym_{k+e}}^{\sym_n}
.$$}
\item{\rm (ii)}{As $\complex\sym_n$-modules, the $(k-1)$-st homology of
$C_{n, k}$ is isomorphic to the $(k-2)$-nd (co)homology of the complement,
$M_{n, k}^\real$, of the $k$-equal real hyperplane arrangement.}
\endproclaim

\demo{Proof}
We first prove (i).
By Theorem \secd.4, it is enough to show that, for $e \leq n-k$, the
restriction of the character $\chi^{\{\eta(k, e), [n-k-e]\}}$ of $W(D_n)$ to
$\sym_n$ corresponds to the representation $$
\big( 
\id_{n-k-e} \otimes \bigwedge^{k-1} E_{k+e}
\big) \uparrow_{\sym_{n-k-e} \times \sym_{k+e}}^{\sym_n}
$$ of $\sym_n$.

By \cite{{\bf 7}, Proposition 5.4.12}, the character of the $\sym_{k+e}$ module
$\bigwedge^{k-1} E_{k+e}$ is given by the partition $\mu = [e+1, 1^{k-1}]$.
The character of the $\sym_{n-k-e}$-module $\id_{n-k-e}$ is given by the
one-row partition $\nu = [n-k-e]$.  Using standard results \cite{{\bf 7}, Definition
6.1.1},
the character of the induction product of these two characters to $\sym_n$
is $$
\sum_\lambda c_{\mu \nu}^\lambda \chi^\lambda
.$$  The proof of (i) is completed by Lemma \secb.1 (i), which shows that 
we also have $$
\chi^{\{\eta(k, e), [n-k-e]\}}
\downarrow^{W(D_n)}_{\sym_n}
= \chi^{\{\mu, \nu\}}
\downarrow^{W(D_n)}_{\sym_n}
= 
\sum_\lambda c_{\mu \nu}^\lambda \chi^\lambda
.$$

Substituting $s = 1$ into \cite{{\bf 14}, Theorem 4.4}, we see that the 
complex character of the $(k-2)$-nd cohomology of $M_{n, k}^{\real}$, regarded
as a $\complex \sym_n$-module, agrees with the character of the representation
described in part (i).  This completes the proof of (ii).
\qed\enddemo

\remark{Remark \secd.8}
Note that, under the usual identifications, we have 
$$W(D_{n-1}) \cap W(A_n) = W(A_{n-1}).$$  It follows that the type $A$ homology
representations described in Theorem \secd.7 have a branching rule analogous
to the type $D$ branching rule of Corollary \secd.6.  However, this would not
be such an obvious result in the absence of the wider context of the type
$D$ representations.
\endremark

\head Acknowledgements \endhead

I thank Markus Pflaum and Nat Thiem for some helpful conversations.

\leftheadtext{} \rightheadtext{}

\end